\documentclass{article}
 
\pdfpagewidth\paperwidth
\pdfpageheight\paperheight
\usepackage{graphicx}
\usepackage{fancyhdr}
\usepackage{hyperref}
\begin{document}
\parindent 0cm
\begin{center}
\baselineskip 14pt

\textbf{\Large{SYMMETRIES RELATED TO DOMINO}}
\vspace{2pt}

\textbf{\Large{TILINGS ON A CHESSBOARD}}
\vspace{24pt}
\baselineskip 16pt

\large{Mih\'aly Hujter $^{1,*}$, Andr\'as Kaszanyitzky $^2$}
\end{center}
\vspace{48pt}
\baselineskip 10pt
\footnotesize{$^1$ \textit{Address}: 
Budapest University of Technology and Economics, Institute of Mathematics, 
H-1521, Budapest, Egry J\'ozsef utca 1.
\\E-mail: hujter@math.bme.hu}

\vspace{2pt}
\footnotesize{$^2$ \textit{Address}: 
Haydn Barytontrio Budapest.\\
\url{http://haydnbarytontrio.hu/en_index.htm}\\
E-mail: kaszi75@gmail.com}
 
\vspace{2pt}
$^*$ \footnotesize{Corresponding author}

\vspace{24pt}

\baselineskip 14pt
\normalsize

\textbf{Abstract}:
\textit{In this paper we study different kinds of symmetries related to 
the domino tilings of chessboards.}

\vspace{12pt}



\baselineskip 14pt
\setlength{\parskip}{12pt}

\vspace{12pt}

\textbf{\large{1. PARTITIONING THE DOMINOS}}

In 1961 Kasteleyn and independently Temperley and Fischer proved that the number of
different domino tilings of a $(2r)\times(2n)$ chessboard is exactly
\[
\prod_{m=1}^r\prod_{k=1}^n
\left(4\cos^2\frac{m\pi}{2r+1}+4\cos^2\frac{k\pi}{2n+1}\right) \ .  
\]
This formula is well-known as the \emph{Kasteleyn} formula. For example, 
for the $6\times6$ chessboard we have $6728$ different tilings.
As an example, we show one of them here. 

\setlength{\unitlength}{0.55mm}
\begin{center}
\begin{picture}(60,60)
\put( 01, 01 ){\line(1,0){ 8}}
\put( 01, 01 ){\line(0,1){18}}
\put( 01, 19 ){\line(1,0){ 8}}
\put( 09, 01 ){\line(0,1){18}}

\put( 11, 01 ){\line(1,0){ 8}}
\put( 11, 01 ){\line(0,1){18}}
\put( 11, 19 ){\line(1,0){ 8}}
\put( 19, 01 ){\line(0,1){18}}

\put( 01, 21 ){\line(1,0){ 8}}
\put( 01, 21 ){\line(0,1){18}}
\put( 01, 39 ){\line(1,0){ 8}}
\put( 09, 21 ){\line(0,1){18}}

\put( 21, 41 ){\line(1,0){ 8}}
\put( 21, 41 ){\line(0,1){18}}
\put( 21, 59 ){\line(1,0){ 8}}
\put( 29, 41 ){\line(0,1){18}}

\put( 31, 41 ){\line(1,0){ 8}}
\put( 31, 41 ){\line(0,1){18}}
\put( 31, 59 ){\line(1,0){ 8}}
\put( 39, 41 ){\line(0,1){18}}

\put( 31, 21 ){\line(1,0){ 8}}
\put( 31, 21 ){\line(0,1){18}}
\put( 31, 39 ){\line(1,0){ 8}}
\put( 39, 21 ){\line(0,1){18}}

\put( 41, 31 ){\line(1,0){ 8}}
\put( 41, 31 ){\line(0,1){18}}
\put( 41, 49 ){\line(1,0){ 8}}
\put( 49, 31 ){\line(0,1){18}}

\put( 51, 31 ){\line(1,0){ 8}}
\put( 51, 31 ){\line(0,1){18}}
\put( 51, 49 ){\line(1,0){ 8}}
\put( 59, 31 ){\line(0,1){18}}

\put( 51, 11 ){\line(1,0){ 8}}
\put( 51, 11 ){\line(0,1){18}}
\put( 51, 29 ){\line(1,0){ 8}}
\put( 59, 11 ){\line(0,1){18}}

\put( 41, 11 ){\line(1,0){ 8}}
\put( 41, 11 ){\line(0,1){18}}
\put( 41, 29 ){\line(1,0){ 8}}
\put( 49, 11 ){\line(0,1){18}}

\put( 11, 21 ){\line(1,0){18}}
\put( 11, 21 ){\line(0,1){ 8}}
\put( 11, 29 ){\line(1,0){18}}
\put( 29, 21 ){\line(0,1){ 8}}

\put( 11, 31 ){\line(1,0){18}}
\put( 11, 31 ){\line(0,1){ 8}}
\put( 11, 39 ){\line(1,0){18}}
\put( 29, 31 ){\line(0,1){ 8}}

\put( 01, 41 ){\line(1,0){18}}
\put( 01, 41 ){\line(0,1){ 8}}
\put( 01, 49 ){\line(1,0){18}}
\put( 19, 41 ){\line(0,1){ 8}}

\put( 01, 51 ){\line(1,0){18}}
\put( 01, 51 ){\line(0,1){ 8}}
\put( 01, 59 ){\line(1,0){18}}
\put( 19, 51 ){\line(0,1){ 8}}

\put( 21, 11 ){\line(1,0){18}}
\put( 21, 11 ){\line(0,1){ 8}}
\put( 21, 19 ){\line(1,0){18}}
\put( 39, 11 ){\line(0,1){ 8}}

\put( 21, 01 ){\line(1,0){18}}
\put( 21, 01 ){\line(0,1){ 8}}
\put( 21, 09 ){\line(1,0){18}}
\put( 39, 01 ){\line(0,1){ 8}}

\put( 41, 01 ){\line(1,0){18}}
\put( 41, 01 ){\line(0,1){ 8}}
\put( 41, 09 ){\line(1,0){18}}
\put( 59, 01 ){\line(0,1){ 8}}

\put( 41, 51 ){\line(1,0){18}}
\put( 41, 51 ){\line(0,1){ 8}}
\put( 41, 59 ){\line(1,0){18}}
\put( 59, 51 ){\line(0,1){ 8}}

\end{picture}
\end{center}

Note that $6728=2^329^2$ is really a nice number; 
however, none of Kasteleyn's eighteen 
$\cos^2$ terms is constructible with ruler and compass.
Interestingly enough, in case of a $16\times16$ chessboard, all terms of the 
Kasteleyn formula is constructible. (The interested reader can find the details 
in our recent paper listed in the references.)

It is a well-known math competition problem that for any domino tiling 
out of the $6728$ different variations, we can always find a straight line 
completely separating the dominos into two nonempty groups. 
(No domino is bisectible, of course.)
Our next figure shows the solution for our example above.    

\begin{center}
\begin{picture}(60,60)
\put( 01, 01 ){\line(1,0){ 8}}
\put( 01, 01 ){\line(0,1){18}}
\put( 01, 19 ){\line(1,0){ 8}}
\put( 09, 01 ){\line(0,1){18}}

\put( 11, 01 ){\line(1,0){ 8}}
\put( 11, 01 ){\line(0,1){18}}
\put( 11, 19 ){\line(1,0){ 8}}
\put( 19, 01 ){\line(0,1){18}}

\put( 01, 21 ){\line(1,0){ 8}}
\put( 01, 21 ){\line(0,1){18}}
\put( 01, 39 ){\line(1,0){ 8}}
\put( 09, 21 ){\line(0,1){18}}

\put( 21, 41 ){\line(1,0){ 8}}
\put( 21, 41 ){\line(0,1){18}}
\put( 21, 59 ){\line(1,0){ 8}}
\put( 29, 41 ){\line(0,1){18}}

\put( 31, 41 ){\line(1,0){ 8}}
\put( 31, 41 ){\line(0,1){18}}
\put( 31, 59 ){\line(1,0){ 8}}
\put( 39, 41 ){\line(0,1){18}}

\put( 31, 21 ){\line(1,0){ 8}}
\put( 31, 21 ){\line(0,1){18}}
\put( 31, 39 ){\line(1,0){ 8}}
\put( 39, 21 ){\line(0,1){18}}

\put( 41, 31 ){\line(1,0){ 8}}
\put( 41, 31 ){\line(0,1){18}}
\put( 41, 49 ){\line(1,0){ 8}}
\put( 49, 31 ){\line(0,1){18}}

\put( 51, 31 ){\line(1,0){ 8}}
\put( 51, 31 ){\line(0,1){18}}
\put( 51, 49 ){\line(1,0){ 8}}
\put( 59, 31 ){\line(0,1){18}}

\put( 51, 11 ){\line(1,0){ 8}}
\put( 51, 11 ){\line(0,1){18}}
\put( 51, 29 ){\line(1,0){ 8}}
\put( 59, 11 ){\line(0,1){18}}

\put( 41, 11 ){\line(1,0){ 8}}
\put( 41, 11 ){\line(0,1){18}}
\put( 41, 29 ){\line(1,0){ 8}}
\put( 49, 11 ){\line(0,1){18}}

\put( 11, 21 ){\line(1,0){18}}
\put( 11, 21 ){\line(0,1){ 8}}
\put( 11, 29 ){\line(1,0){18}}
\put( 29, 21 ){\line(0,1){ 8}}

\put( 11, 31 ){\line(1,0){18}}
\put( 11, 31 ){\line(0,1){ 8}}
\put( 11, 39 ){\line(1,0){18}}
\put( 29, 31 ){\line(0,1){ 8}}

\put( 01, 41 ){\line(1,0){18}}
\put( 01, 41 ){\line(0,1){ 8}}
\put( 01, 49 ){\line(1,0){18}}
\put( 19, 41 ){\line(0,1){ 8}}

\put( 01, 51 ){\line(1,0){18}}
\put( 01, 51 ){\line(0,1){ 8}}
\put( 01, 59 ){\line(1,0){18}}
\put( 19, 51 ){\line(0,1){ 8}}

\put( 21, 11 ){\line(1,0){18}}
\put( 21, 11 ){\line(0,1){ 8}}
\put( 21, 19 ){\line(1,0){18}}
\put( 39, 11 ){\line(0,1){ 8}}

\put( 21, 01 ){\line(1,0){18}}
\put( 21, 01 ){\line(0,1){ 8}}
\put( 21, 09 ){\line(1,0){18}}
\put( 39, 01 ){\line(0,1){ 8}}

\put( 41, 01 ){\line(1,0){18}}
\put( 41, 01 ){\line(0,1){ 8}}
\put( 41, 09 ){\line(1,0){18}}
\put( 59, 01 ){\line(0,1){ 8}}

\put( 41, 51 ){\line(1,0){18}}
\put( 41, 51 ){\line(0,1){ 8}}
\put( 41, 59 ){\line(1,0){18}}
\put( 59, 51 ){\line(0,1){ 8}}

\put( 40, 00 ){\line(0,1){60}}

\end{picture}
\end{center}

In the recent paper we ask a similar question. 
As in graph theory, by a Hamiltonian path, we mean 
such a path which starts at a vertex, ends at another vertex, and
goes through all other vertices. 
In summary, a Hamiltonian path visits each vertex exactly once. 
(Here we consider each vertex of any elementary square of the chessboard.)
On a $6\times6$ chessboard the number of vertices is $49$. 
Therefore, if there exists any Hamiltonian path, the length of each Hamiltonian 
path is $48$. Our next figure shows a Hamiltonian path connecting two opposite 
corners of the chessboard.  
 
\begin{center}
\begin{picture}(140,60)

\put( 01, 01 ){\line(1,0){ 8}}
\put( 01, 01 ){\line(0,1){18}}
\put( 01, 19 ){\line(1,0){ 8}}
\put( 09, 01 ){\line(0,1){18}}

\put( 11, 01 ){\line(1,0){ 8}}
\put( 11, 01 ){\line(0,1){18}}
\put( 11, 19 ){\line(1,0){ 8}}
\put( 19, 01 ){\line(0,1){18}}

\put( 01, 21 ){\line(1,0){ 8}}
\put( 01, 21 ){\line(0,1){18}}
\put( 01, 39 ){\line(1,0){ 8}}
\put( 09, 21 ){\line(0,1){18}}

\put( 21, 41 ){\line(1,0){ 8}}
\put( 21, 41 ){\line(0,1){18}}
\put( 21, 59 ){\line(1,0){ 8}}
\put( 29, 41 ){\line(0,1){18}}

\put( 31, 41 ){\line(1,0){ 8}}
\put( 31, 41 ){\line(0,1){18}}
\put( 31, 59 ){\line(1,0){ 8}}
\put( 39, 41 ){\line(0,1){18}}

\put( 31, 21 ){\line(1,0){ 8}}
\put( 31, 21 ){\line(0,1){18}}
\put( 31, 39 ){\line(1,0){ 8}}
\put( 39, 21 ){\line(0,1){18}}

\put( 41, 31 ){\line(1,0){ 8}}
\put( 41, 31 ){\line(0,1){18}}
\put( 41, 49 ){\line(1,0){ 8}}
\put( 49, 31 ){\line(0,1){18}}

\put( 51, 31 ){\line(1,0){ 8}}
\put( 51, 31 ){\line(0,1){18}}
\put( 51, 49 ){\line(1,0){ 8}}
\put( 59, 31 ){\line(0,1){18}}

\put( 51, 11 ){\line(1,0){ 8}}
\put( 51, 11 ){\line(0,1){18}}
\put( 51, 29 ){\line(1,0){ 8}}
\put( 59, 11 ){\line(0,1){18}}

\put( 41, 11 ){\line(1,0){ 8}}
\put( 41, 11 ){\line(0,1){18}}
\put( 41, 29 ){\line(1,0){ 8}}
\put( 49, 11 ){\line(0,1){18}}

\put( 11, 21 ){\line(1,0){18}}
\put( 11, 21 ){\line(0,1){ 8}}
\put( 11, 29 ){\line(1,0){18}}
\put( 29, 21 ){\line(0,1){ 8}}

\put( 11, 31 ){\line(1,0){18}}
\put( 11, 31 ){\line(0,1){ 8}}
\put( 11, 39 ){\line(1,0){18}}
\put( 29, 31 ){\line(0,1){ 8}}

\put( 01, 41 ){\line(1,0){18}}
\put( 01, 41 ){\line(0,1){ 8}}
\put( 01, 49 ){\line(1,0){18}}
\put( 19, 41 ){\line(0,1){ 8}}

\put( 01, 51 ){\line(1,0){18}}
\put( 01, 51 ){\line(0,1){ 8}}
\put( 01, 59 ){\line(1,0){18}}
\put( 19, 51 ){\line(0,1){ 8}}

\put( 21, 11 ){\line(1,0){18}}
\put( 21, 11 ){\line(0,1){ 8}}
\put( 21, 19 ){\line(1,0){18}}
\put( 39, 11 ){\line(0,1){ 8}}

\put( 21, 01 ){\line(1,0){18}}
\put( 21, 01 ){\line(0,1){ 8}}
\put( 21, 09 ){\line(1,0){18}}
\put( 39, 01 ){\line(0,1){ 8}}

\put( 41, 01 ){\line(1,0){18}}
\put( 41, 01 ){\line(0,1){ 8}}
\put( 41, 09 ){\line(1,0){18}}
\put( 59, 01 ){\line(0,1){ 8}}

\put( 41, 51 ){\line(1,0){18}}
\put( 41, 51 ){\line(0,1){ 8}}
\put( 41, 59 ){\line(1,0){18}}
\put( 59, 51 ){\line(0,1){ 8}}

\put( 00, 00 ){\line(0,1){40}}
\put( 10, 00 ){\line(0,1){30}}
\put( 50, 10 ){\line(0,1){40}}
\put( 60, 00 ){\line(0,1){50}}
\put( 30, 30 ){\line(0,1){30}}
\put( 00, 50 ){\line(0,1){10}}
\put( 20, 40 ){\line(0,1){10}}
\put( 40, 20 ){\line(0,1){40}}
\put( 20, 10 ){\line(0,1){10}}
\put( 10, 00 ){\line(1,0){50}}
\put( 20, 10 ){\line(1,0){30}}
\put( 10, 30 ){\line(1,0){20}}
\put( 00, 50 ){\line(1,0){20}}
\put( 00, 60 ){\line(1,0){30}}
\put( 40, 60 ){\line(1,0){20}}
\put( 50, 50 ){\line(1,0){10}}
\put( 20, 20 ){\line(1,0){20}}
\put( 00, 40 ){\line(1,0){20}}

\put( 80, 00 ){\line(0,1){40}}
\put( 90, 00 ){\line(0,1){30}}
\put(130, 10 ){\line(0,1){40}}
\put(140, 00 ){\line(0,1){50}}
\put(110, 30 ){\line(0,1){30}}
\put( 80, 50 ){\line(0,1){10}}
\put(100, 40 ){\line(0,1){10}}
\put(120, 20 ){\line(0,1){40}}
\put(100, 10 ){\line(0,1){10}}
\put( 90, 00 ){\line(1,0){50}}
\put(100, 10 ){\line(1,0){30}}
\put( 90, 30 ){\line(1,0){20}}
\put( 80, 50 ){\line(1,0){20}}
\put( 80, 60 ){\line(1,0){30}}
\put(120, 60 ){\line(1,0){20}}
\put(130, 50 ){\line(1,0){10}}
\put(100, 20 ){\line(1,0){20}}
\put( 80, 40 ){\line(1,0){20}}
\put( 80, 00 ){\line(0,1){40}}
\put( 90, 00 ){\line(0,1){30}}
\put(130, 10 ){\line(0,1){40}}
\put(140, 00 ){\line(0,1){50}}
\put(110, 30 ){\line(0,1){30}}
\put( 80, 50 ){\line(0,1){10}}
\put(100, 40 ){\line(0,1){10}}
\put(120, 20 ){\line(0,1){40}}
\put(100, 10 ){\line(0,1){10}}
\put( 90, 00 ){\line(1,0){50}}
\put(100, 10 ){\line(1,0){30}}
\put( 90, 30 ){\line(1,0){20}}
\put( 80, 50 ){\line(1,0){20}}
\put( 80, 60 ){\line(1,0){30}}
\put(120, 60 ){\line(1,0){20}}
\put(130, 50 ){\line(1,0){10}}
\put(100, 20 ){\line(1,0){20}}
\put( 80, 40 ){\line(1,0){20}}

\put( 80, 00  ){\circle*{2}}
\put( 80, 10  ){\circle*{2}}
\put( 80, 20  ){\circle*{2}}
\put( 80, 30  ){\circle*{2}}
\put( 80, 40  ){\circle*{2}}
\put( 80, 50  ){\circle*{2}}
\put( 80, 60  ){\circle*{2}}

\put( 90, 00  ){\circle*{2}}
\put( 90, 10  ){\circle*{2}}
\put( 90, 20  ){\circle*{2}}
\put( 90, 30  ){\circle*{2}}
\put( 90, 40  ){\circle*{2}}
\put( 90, 50  ){\circle*{2}}
\put( 90, 60  ){\circle*{2}}

\put(100, 00  ){\circle*{2}}
\put(100, 10  ){\circle*{2}}
\put(100, 20  ){\circle*{2}}
\put(100, 30  ){\circle*{2}}
\put(100, 40  ){\circle*{2}}
\put(100, 50  ){\circle*{2}}
\put(100, 60  ){\circle*{2}}

\put(110, 00  ){\circle*{2}}
\put(110, 10  ){\circle*{2}}
\put(110, 20  ){\circle*{2}}
\put(110, 30  ){\circle*{2}}
\put(110, 40  ){\circle*{2}}
\put(110, 50  ){\circle*{2}}
\put(110, 60  ){\circle*{2}}

\put(120, 00  ){\circle*{2}}
\put(120, 10  ){\circle*{2}}
\put(120, 20  ){\circle*{2}}
\put(120, 30  ){\circle*{2}}
\put(120, 40  ){\circle*{2}}
\put(120, 50  ){\circle*{2}}
\put(120, 60  ){\circle*{2}}

\put(130, 00  ){\circle*{2}}
\put(130, 10  ){\circle*{2}}
\put(130, 20  ){\circle*{2}}
\put(130, 30  ){\circle*{2}}
\put(130, 40  ){\circle*{2}}
\put(130, 50  ){\circle*{2}}
\put(130, 60  ){\circle*{2}}

\put(140, 00  ){\circle*{2}}
\put(140, 10  ){\circle*{2}}
\put(140, 20  ){\circle*{2}}
\put(140, 30  ){\circle*{2}}
\put(140, 40  ){\circle*{2}}
\put(140, 50  ){\circle*{2}}
\put(140, 60  ){\circle*{2}}

\end{picture}
\end{center}

We give the reader an exercise: 
Find a Hamiltonian path which connects
the top-left and the bottom-right corners.

Or next figure gives a larger chessboard tiling with 
a corresponding Hamiltonian path. On both sides of the 
Hamiltonian path we find $16$ dominos. 

\begin{figure}[ht]
\centerline{
\includegraphics[width=0.55\textwidth]{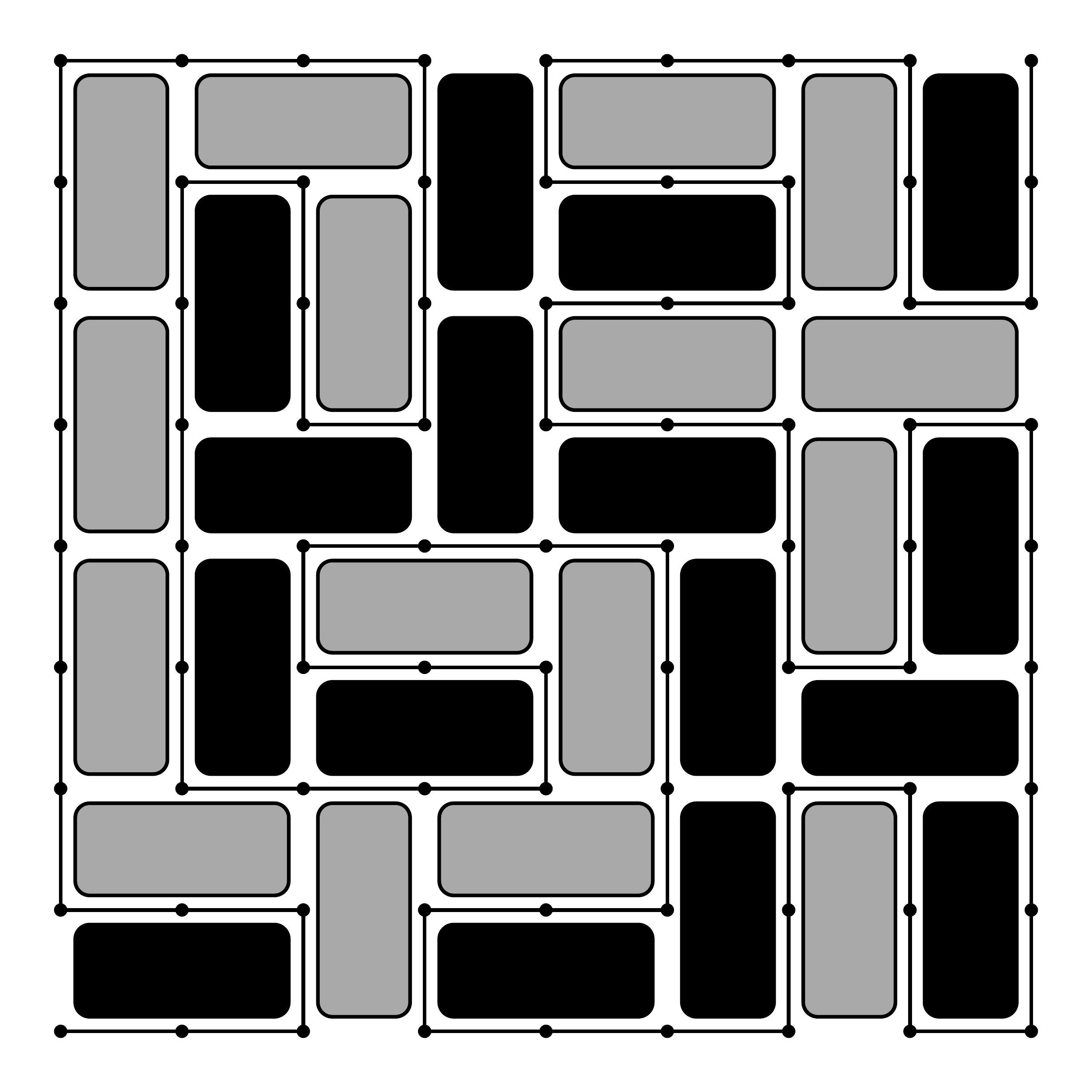}
}
\end{figure} 

The chessboard can be any rectangle with sides of even lengths.
The next figures show a chessboard with two Hamiltonian paths. 
The domino tilings are the same in both cases. 
Note that the Hamiltonian path connecting the first and the third 
corners is very different from the Hamiltonian path connecting 
the second and the fourth corners of the chessboard. 

Given an arbitrary domino tiling, 
in our recent paper listed in the references we showed 
the existence of a Hamiltonian path between the bottom-left 
and the top-right corners. 
Similarly, there is another  
Hamiltonian path between the top-left and the bottom-right 
corners. 
Furthermore, these two Hamiltonian paths both have 
the following nice property. 
Assume that we start from either the bottom-left 
or the bottom-right corner. 
Then both Hamiltonian paths keep the one-way traffic rule that the
first (i.e.\ the left-most), the third, the fifth, etc.\ 
vertical lines all go to the north, 
and the second, the fourth, the sixth, etc.\ vertical 
lines all go to the south. 
Considering the horizontal lines, 
the first Hamiltonian path goes to the east on the bottom line, 
on the third line, on the fifth line, etc., and goes to west 
on the second line, on the fourth line, 
etc.; on the other hand, the second Hamiltonian path    
goes to the west on the bottom line, on the third line, 
on the fifth line, etc., 
and goes to east on the second line, on the fourth line, 
and so on. 

\begin{figure}[ht]
\centerline{
\includegraphics[width=0.5\textwidth]{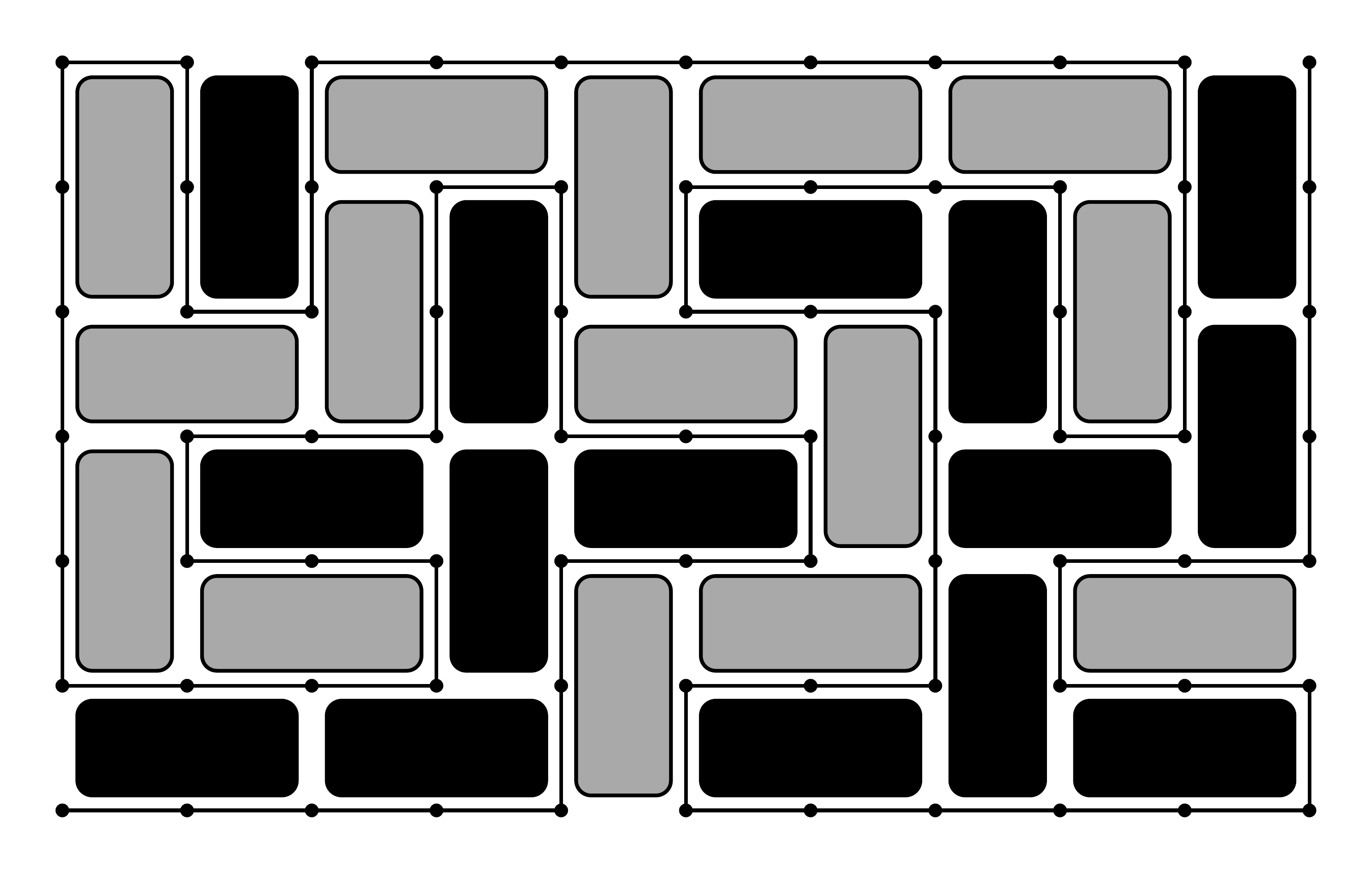}
\ \ 
\includegraphics[width=0.5\textwidth]{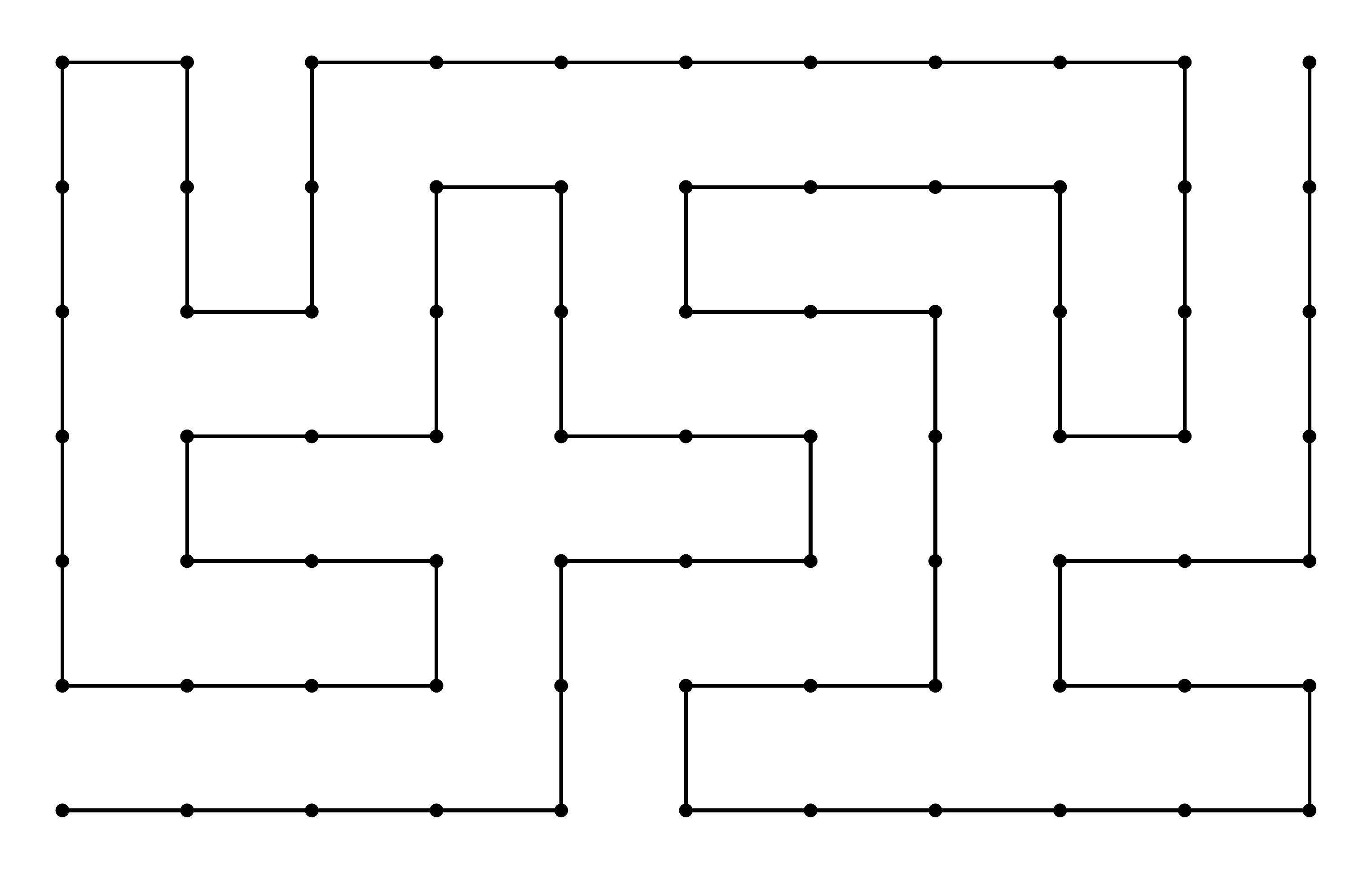}
}
\smallskip
\centerline{
\includegraphics[width=0.5\textwidth]{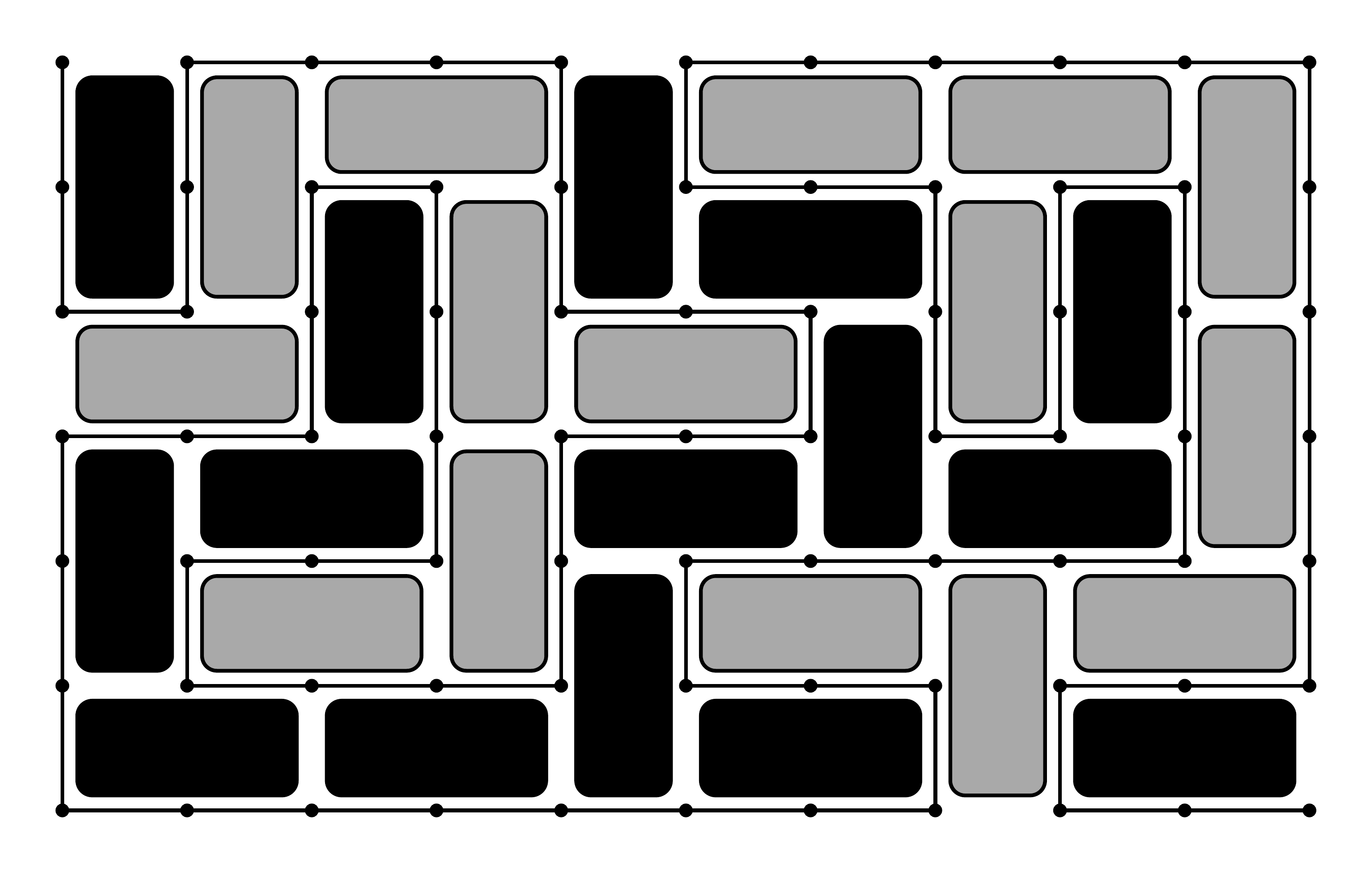}
\ \  
\includegraphics[width=0.5\textwidth]{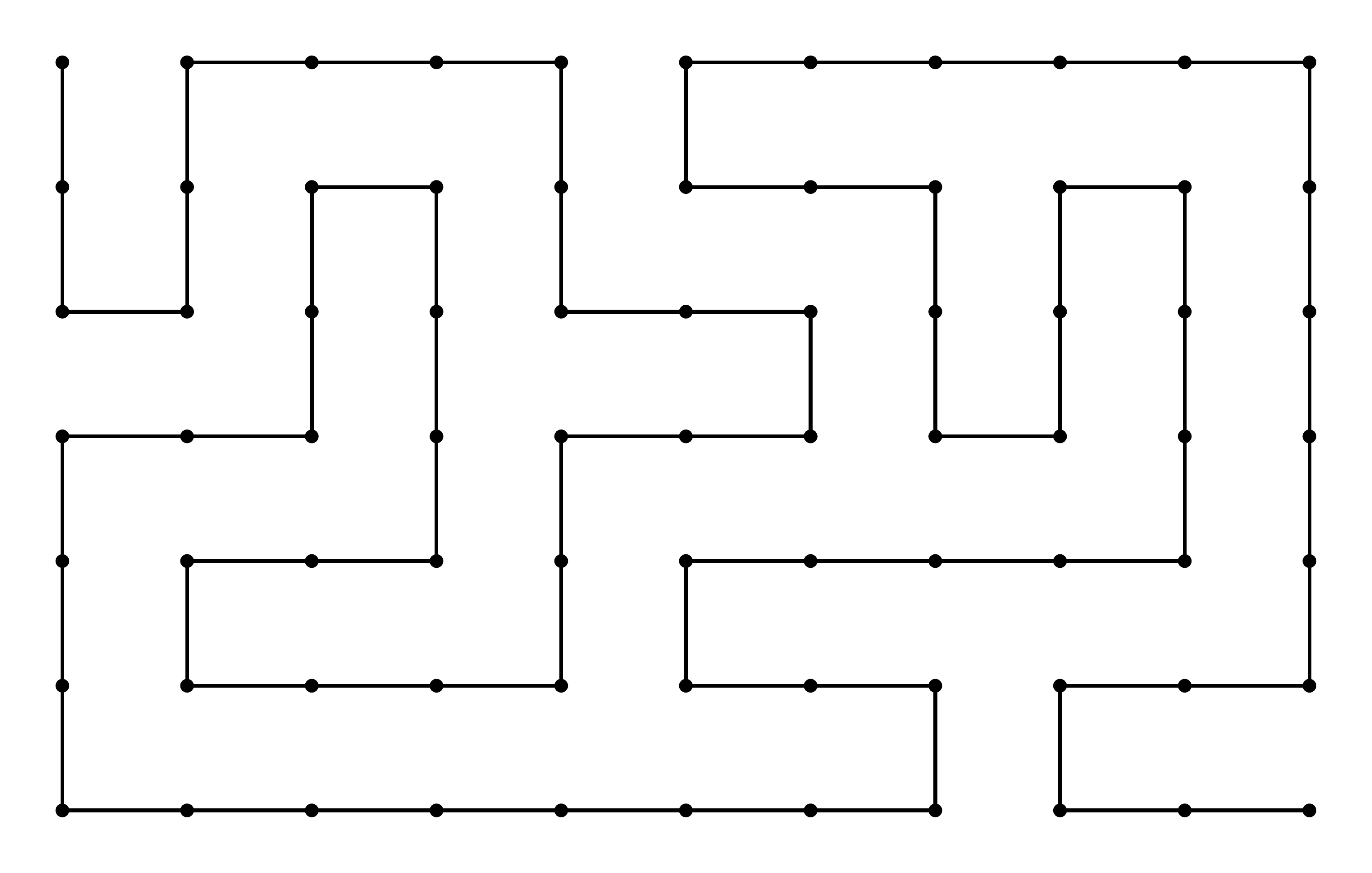}
}
\end{figure}

\newpage

\textbf{\large{2. GENERATING FUNCTION AND NEW RECURRENCE RELATION}}

Let $c_n$ denote the number of domino tilings of a $6\times(2n)$
chessboard. 

As a recent result, from the OEIS Foundation Inc., we learn 
the generating function of the sequence $c_0,c_1,c_2,\ldots$:
\begin{eqnarray*}
g_{6}(x) &=&\frac{1-27x+177x^2-328x^3+177x^4-27x^5+x^6}
{(1-x)\left(1-39x+377x^2-847x^3+377x^4-39x^5+x^6\right)} \\
&=&1+13x+281x^2+6728x^3+167089x^4+4213133x^5+...
\end{eqnarray*}%
Observe that 
\[
(1-x)g_6(x)=\frac{u^3-27u^2+174u-274}{u^3-39u^2+374u-769} 
\]
where $u=x+x^{-1}$.

By using the above results we have the following values.
\[
\begin{tabular}{|c|cccccccc|}
\hline
$n$ & $0$ & $1$ & $2$ & $3$ & $4$ & $5$ & $6$ & $7$ \\ \hline
$c_n$ & $1$ & $13$ & $281$ & $6728$ & $167089$ & $4213133$ & $106912793$ & 
$2720246633$ \\ \hline
\end{tabular}
\]

In some sense this method of computing the numbers $c_0,c_1,c_2,\ldots$
is superior to the Kasteleyn formula because you do not have 
to worry about the numerical errors.
\[
\begin{tabular}{|c|ccc|}
\hline
$n$ & $8$ & $9$ & $10$ \\ \hline
$c_n$ & $69289288909$ & $1765722581057$ & $45005025662792$ \\ \hline
\end{tabular}
\]
\[
\begin{tabular}{|c|ccc|}
\hline
$n$ & $11$ & $12$ & $13$ \\ \hline
$c_n$ & $\ 1147185247901449$ & $29242880940226381$ & $745439797095329713$
\\ \hline
\end{tabular}
\]
\[
\begin{tabular}{|c|cc|}
\hline
$n$ & $14$ & $15$ \\ \hline
$c_n$ & $19002353776441540177$ & $484398978524471931341$ \\ \hline
\end{tabular}
\]
\[
\begin{tabular}{|c|cc|}
\hline
$n$ & $16$ & $17$ \\ \hline
$c_n$ & $12348080425980866090537$ & $314771823879840325570888$ \\ \hline
\end{tabular}
\]

In general, we can get these values as the top-left entry of the matrix 
$C^n$ where $C$ is the following $20\times 20$ symmetric matrix.

\[
C=\left[ 
\begin{array}{cccccccccccccccccccc}
13 & 5 & 3 & 4 & 3 & 5 & 2 & 1 & 2 & 1 & 2 & 1 & 2 & 1 & 2 & 1 & 1 & 0 & 1 & 
1 \\ 
5 & 5 & 0 & 2 & 1 & 2 & 2 & 0 & 1 & 1 & 2 & 0 & 1 & 1 & 2 & 1 & 0 & 0 & 1 & 1
\\ 
3 & 0 & 3 & 0 & 1 & 1 & 0 & 1 & 0 & 0 & 0 & 1 & 0 & 0 & 0 & 0 & 1 & 0 & 0 & 0
\\ 
4 & 2 & 0 & 4 & 0 & 2 & 2 & 0 & 2 & 1 & 2 & 0 & 2 & 1 & 1 & 1 & 0 & 0 & 0 & 1
\\ 
3 & 1 & 1 & 0 & 3 & 0 & 0 & 1 & 0 & 0 & 0 & 1 & 0 & 0 & 0 & 0 & 0 & 0 & 1 & 0
\\ 
5 & 2 & 1 & 2 & 0 & 5 & 1 & 0 & 2 & 1 & 1 & 0 & 2 & 1 & 2 & 1 & 1 & 0 & 0 & 1
\\ 
2 & 2 & 0 & 2 & 0 & 1 & 4 & 0 & 1 & 2 & 2 & 0 & 1 & 2 & 1 & 1 & 0 & 0 & 0 & 1
\\ 
1 & 0 & 1 & 0 & 1 & 0 & 0 & 2 & 0 & 0 & 0 & 1 & 0 & 0 & 0 & 0 & 0 & 0 & 0 & 0
\\ 
2 & 1 & 0 & 2 & 0 & 2 & 1 & 0 & 4 & 2 & 1 & 0 & 2 & 1 & 1 & 2 & 0 & 0 & 0 & 1
\\ 
1 & 1 & 0 & 1 & 0 & 1 & 2 & 0 & 2 & 5 & 1 & 0 & 1 & 2 & 1 & 2 & 0 & 1 & 0 & 1
\\ 
2 & 2 & 0 & 2 & 0 & 1 & 2 & 0 & 1 & 1 & 2 & 0 & 1 & 1 & 1 & 1 & 0 & 0 & 0 & 1
\\ 
1 & 0 & 1 & 0 & 1 & 0 & 0 & 1 & 0 & 0 & 0 & 1 & 0 & 0 & 0 & 0 & 0 & 0 & 0 & 0
\\ 
2 & 1 & 0 & 2 & 0 & 2 & 1 & 0 & 2 & 1 & 1 & 0 & 2 & 1 & 1 & 1 & 0 & 0 & 0 & 1
\\ 
1 & 1 & 0 & 1 & 0 & 1 & 2 & 0 & 1 & 2 & 1 & 0 & 1 & 2 & 1 & 1 & 0 & 0 & 0 & 1
\\ 
2 & 2 & 0 & 1 & 0 & 2 & 1 & 0 & 1 & 1 & 1 & 0 & 1 & 1 & 2 & 1 & 0 & 0 & 0 & 1
\\ 
1 & 1 & 0 & 1 & 0 & 1 & 1 & 0 & 2 & 2 & 1 & 0 & 1 & 1 & 1 & 2 & 0 & 0 & 0 & 1
\\ 
1 & 0 & 1 & 0 & 0 & 1 & 0 & 0 & 0 & 0 & 0 & 0 & 0 & 0 & 0 & 0 & 1 & 0 & 0 & 0
\\ 
0 & 0 & 0 & 0 & 0 & 0 & 0 & 0 & 0 & 1 & 0 & 0 & 0 & 0 & 0 & 0 & 0 & 1 & 0 & 0
\\ 
1 & 1 & 0 & 0 & 1 & 0 & 0 & 0 & 0 & 0 & 0 & 0 & 0 & 0 & 0 & 0 & 0 & 0 & 1 & 0
\\ 
1 & 1 & 0 & 1 & 0 & 1 & 1 & 0 & 1 & 1 & 1 & 0 & 1 & 1 & 1 & 1 & 0 & 0 & 0 & 1
\end{array}
\right] 
\]

We call our brand new method {\em the matrix power method for computing
the Kasteleyn formula}. The proof of the correctness of this method
is going to be published in a forthcoming paper.
Note that the matrix power method enables us to compute the 
4096th element, for example, much faster than the 
Kasteleyn formula or the generating function method.  
First you compute $C^2$, then $C^4$ as the square of $C^2$,
then $C^8$ as the square of $C^4$, and so on.

Interestingly enough, the characteristic polynomial of both $C$ and its
inverse matrix $C^{-1}$ can be written in the following form.

\begin{eqnarray*}
&&\lambda ^{10}(\left(\lambda^{10}+\lambda^{-10}\right)  \\
&&-63(\lambda^9+\lambda^{-9})+1561(\lambda^8+\lambda^{-8})
-21023\left( \lambda^7+\lambda ^{-7}\right)  \\
&&+176393\allowbreak\left(\lambda^6+\lambda ^{-6}\right)
-992383(\lambda^5+\lambda^{-5})+3912609(\lambda^4+\lambda^{-4}) \\
&&-11117\,602(\lambda^3+\lambda^{-3})+23182782
\left(\lambda^2+\lambda^{-2}\right)-35879970(\lambda+\lambda ^{-1}) \\
&&+41475390)
\end{eqnarray*}

A corollary of this fact is that the numbers $c_n$ satisfy the following
recurrence relation (for $n\geq 10$):
\begin{eqnarray*}
&&c_{n+10}=-c_{n-10} \\
&&+63(c_{n-9}+c_{n+9})-1561(c_{n-8}+c_{n+8})+21023\left(
c_{n-7}+c_{n+7}\right)  \\
&&-176393\left(c_{n-6}+c_{n+6}\right)
+992383(c_{n-5}+c_{n+5})-3912609(c_{n-4}+c_{n+4}) \\
&&+11117\,602(c_{n-3}+c_{n+3})-23182782\left(c_{n-2}+c_{n+2}\right)
+35879970(c_{n-1}+c_n) \\
&&-41475390c_n
\end{eqnarray*}
This enables us to compute $c_{n+10}$ if the previous twenty numbers 
$c_{n-10},c_{n-9},\allowbreak\ldots,\allowbreak c_n,\allowbreak\ldots,c_{n+9}$ 
are already known.  
Observe that this formula is symmetric with respect to $c_{n-k}$ and $c_{n+k}$
where $k=1,2,\ldots,10$. 

\bigskip

\vspace{12pt}
\textbf{\large{3. CONCLUSION}}

Each of the exponentially many different tilings of dominos 
show some symmetries. 
Here we gave some examples. 
The physical and chemical connections might also be interesting.

\newpage

\vspace{12pt}
\textbf{\large{Acknowledgements}}

We express our thanks to the inspiring 
community of the Mathematics Museum in Budapest for
the interesting discussions of many different subjects.

\vspace{12pt}
\textbf{\large{REFERENCES}}
\vspace{12pt}

\bigskip

Hujter, M.\ and Kaszanyitzky, A. (2015)
Hamiltonian paths on directed grids,
\url{http://arxiv.org/abs/1512.00718}

Hujter, M.\ and Kaszanyitzky, A. (2016)
The matrix power method for computing
the Kasteleyn formula,
\textit{in preparation.}

Kasteleyn, P.W. (1961) 
The statistics of dimers on a lattice, 
\textit{Physica}, 27, 1209--1225.

OEIS Foundation Inc (2016) 
\textit{The on-line Encyclopedia of integer sequences}, 
\newline\url{http://oeis.org}

Temperley, H.N.V.\ and Fisher, M.E. (1961)
Dimer problem in statistical mechanics--an exact result, 
\textit{Phil.\ Mag.}, 68, 1061--1063.

\end{document}